Т.Г. Эргашев

## ОБОБЩЕННЫЕ РЕШЕНИЯ ОДНОГО ВЫРОЖДАЮЩЕГОСЯ ГИПЕРБОЛИЧЕСКОГО УРАВНЕНИЯ ВТОРОГО РОДА СО СПЕКТРАЛЬНЫМ ПАРАМЕТРОМ

Для вырождающегося гиперболического уравнения второго рода со спектральным параметром изучены задачи Коши, Коши – Гурса и Гурса в новом классе обобщенных решений и приведен пример, показывающий важность введения понятия такого класса. Вводятся в рассмотрение операторы с функциями Бесселя в ядрах. Установленные здесь тождества для этих операторов играли важную роль при получении явных интегральных представлений исследуемых задач.

**Ключевые слова:** *вырождающееся уравнение гиперболического типа второго рода, спектральный параметр, обобщенное решение, оператор с функциями Бесселя в ядрах.*

### 1. Постановка задачи

Рассмотрим уравнение

$$y^m U_{xx} - U_{yy} + \lambda^2 y^m U = 0, \qquad (1)$$

где $m$ – действительное число, причем $-1 < m < 0$, а $\lambda$ – действительное или чисто мнимое постоянное число.

Пусть $a$, $b$, $M$, $N$ – некоторые действительные числа, причем $M < N < +\infty$, $a, b \in [M, N]$; пусть $D$ – область, ограниченная отрезком $AB = \{(x,y): y = 0, a < x < b\}$ и характеристиками $AB: \xi = a$, $BC: \eta = b$ уравнения (1) при $y > 0$. Здесь

$$\xi = x - \frac{2}{m+2} y^{(m+2)/2}, \quad \eta = x + \frac{2}{m+2} y^{(m+2)/2}. \qquad (2)$$

Для уравнения (1) прямая параболического вырождения является особой характеристикой – огибающей обоих семейств характеристик. В зависимости от степени вырождения $m$ предельные значения $\tau(x) \equiv U(x,+0)$, $\nu(x) \equiv U'_y(x,+0)$ могут иметь особенности. Чтобы обеспечить необходимую гладкость решения $U(x,y)$ вне линии характеристического вырождения, необходимо требовать повышенную гладкость функций $\tau(x)$ и $\nu(x)$. С целью ослабить это требование в [1] дано определение и изучены свойства так называемого класса $R_{2k}^{\lambda}$ (здесь и далее $k$ принимает значения $a$ и $b$) обобщенных решений уравнения (1) в области $D$, который при $\lambda = 0$ и $k = 0$ совпадает с классом $R_2$, введенным и изученным И.Л. Каролем [2]. Кроме того, на основе известной формулы классического решения задачи Коши [3] для уравнения (1) в [1] получен явный и удобный для даль-нейших исследований вид обобщенного решения этой же задачи в классе $R_{2k}^{\lambda}$ и исследованы обобщенные решения, для которых $\tau'(x), \nu(x) \in C(a,b)$ вместо требуемого $C^2[a,b]$.

В настоящей работе исследуются задачи Коши – Гурса и Гурса для уравнения (1) в классе обобщенных решений $R_{2k}^{\lambda}$ и приводится пример, показывающий важность введения понятия такого класса. Прежде чем перейти к решению поставленных задач вводятся в рассмотрение некоторые операторы с функциями Бесселя в ядрах. Именно выявленные здесь новые свойства этих операторов играют важную роль при получении явных интегральных представлений исследуемых задач.

**Задача Коши – Гурса.** Требуется найти в области $D$ решение $U(\xi,\eta)$ уравнения (1) из класса $R_{2k}^{\lambda}$, удовлетворяющее начальному условию

$$U'_y(x,+0) = \nu(x) \tag{3}$$

и одному из условий

$$U|_{AB} = \psi_a(x), \ a \leq x \leq (a+b)/2; \tag{4}$$

$$U|_{BC} = \psi_b(x), \ (a+b)/2 \leq x \leq b, \tag{5}$$

где $\nu(x), \psi_a(x), \psi_b(x)$ – заданные функции.

**Задача Гурса.** Требуется найти в области $D$ решение $U(\xi,\eta)$ уравнения (1) из класса $R_{2k}^{\lambda}$, удовлетворяющее условиям (4), (5) и условию согласования

$$\psi_a[(a+b)/2] = \psi_b[(a+b)/2]. \tag{6}$$

В характеристических координатах (2) уравнение (1) переходит в уравнение

$$u_{\xi\eta} - \frac{\beta}{\eta-\xi}(u_\eta - u_\xi) + \frac{1}{4}\lambda^2 u = 0, \tag{7}$$

область $D$ преобразуется в треугольник $\Delta$, ограниченный прямыми $\xi = a, \ \eta = b$ и $\eta = \xi$, а условия (3) – (6) соответственно принимают вид

$$\left(\frac{m+2}{4}\right)^{2\beta}\lim_{\eta\to\xi}(\eta-\xi)^{2\beta}\left(u_\xi - u_\eta\right) = \nu(\xi), \ a < \xi < b; \tag{8}$$

$$u|_{\xi=a} = \varphi_a(\eta), \ a \leq \eta \leq b; \tag{9}$$

$$u|_{\eta=b} = \varphi_b(\xi), \ a \leq \xi \leq b; \tag{10}$$

$$\varphi_a(b) = \varphi_b(a), \tag{11}$$

где

$$\beta = \frac{m}{2(m+2)}, \ -1 < 2\beta < 0; \ u(\xi,\eta) = U\left[\frac{\eta+\xi}{2}, -\left(\frac{\eta-\xi}{2(1-2\beta)}\right)^{1-2\beta}\right];$$

$$\varphi_a(\eta) = \psi_a\left(\frac{a+\eta}{2}\right), \ \varphi_b(\xi) = \psi_b\left(\frac{\xi+b}{2}\right), \ \varphi_k \in C[a,b]\bigcap C^{(1,\alpha_k)}(a,b), \ \alpha_k > -2\beta.$$

### 2. Некоторые операторы с функциями Бесселя в ядрах и их свойства

Введем в рассмотрение операторы

$$A_{abx}^{n,\lambda}[f(x)] \equiv f(x) - \int_a^x f(t)\left(\frac{b-t}{b-x}\right)^n \frac{\partial}{\partial t}J_0\left[\lambda\sqrt{(x-b)(x-t)}\right]dt; \qquad (12)$$

$$B_{abx}^{n,\lambda}[f(x)] \equiv f(x) + \int_a^x f(t)\left(\frac{b-x}{b-t}\right)^{1-n} \frac{\partial}{\partial x}J_0\left[\lambda\sqrt{(b-t)(x-t)}\right]dt, \qquad (13)$$

где $n$ – неотрицательное целое число, а $J_\alpha(z)$ – функция Бесселя первого рода порядка $\alpha$: $J_\alpha(z) = \sum_{k=0}^{\infty} \frac{(-1)^k (z/2)^{2k+\alpha}}{k!\Gamma(k+\alpha+1)}$. Здесь следует особо отметить, что $a \le b$.

При рассмотрении операторов (12) и (13) предположим, что $f(x) \in C(M,N) \cap L_1[M,N]$. При таких предположениях относительно функции $f(x)$ выражения $A_{abx}^{n,\lambda}[f(x)]$ и $B_{abx}^{n,\lambda}[f(x)]$ будут определены в $(M,N)$ и принадлежат классу $C(M,N)$.

Отметим, что операторы (12) и (13) при $a=b$, $n=0$ и $n=1$ введены и исследованы в [4].

Следующая теорема выражает основное свойство этих операторов.

**Теорема.** Если $f(x) \in C[M,N]$, то для любых $a,b \in [M,N]$ и $x \in (M,N)$ справедливы следующие равенства:

$$B_{abx}^{n,\lambda}\left\{A_{abx}^{n,\lambda}[f(x)]\right\} = f(x), \ A_{abx}^{n,\lambda}\left\{B_{abx}^{n,\lambda}[f(x)]\right\} = f(x),$$

т.е. в классе непрерывных на $[M,N]$ функций операторы (12) и (13) являются взаимно обратными.

*Доказательство.* Подействовав на оператор $A_{abx}^{n,\lambda}[f(x)]$ оператором $B_{abx}^{n,\lambda}$, получаем

$$B_{abx}^{n,\lambda}\left\{A_{abx}^{n,\lambda}[f(x)]\right\} = f(x) - (b-x)^{1-n}\int_0^x f(t)L(x,t;\lambda)(b-t)^n dt,$$

где

$$L(x,t;\lambda) = \frac{1}{b-x}\frac{\partial}{\partial t}J_0\left[\lambda\sqrt{(x-b)(x-t)}\right] - (b-t)\frac{\partial}{\partial x}J_0\left[\lambda\sqrt{(b-t)(x-t)}\right] +$$

$$+ \int_t^x \frac{1}{b-s}\frac{\partial}{\partial x}J_0\left[\lambda\sqrt{(b-s)(x-s)}\right]\frac{\partial}{\partial t}J_0\left[\lambda\sqrt{(s-b)(s-t)}\right]ds.$$

После несложных преобразований выражению $L(x,t;\lambda)$ можно придать вид

$$L(x,t;\lambda) = \sum_{k=0}^{\infty}\left(\frac{\lambda}{2}\right)^{2k+2}\frac{(x-t)^k(x-b)^k}{k!}F_k(z),$$

где
$$F_k(z) = \sum_{i=0}^{k+1}\frac{(-1)^i}{i!(k-i+1)!}F(-k,k-i+1;k+1;z),$$
$$k = 0,1,2,..., \ z = (x-t)/(x-b).$$

Здесь $F(a,b;c;z) = \sum_{n=0}^{\infty} \frac{(a)_n (b)_n}{n!(c)_n} z^n$ – гипергеометрическая функция Гаусса [5], а $(a)_n$ – символ Похгаммера [5]:

$$(a)_0 = 1, \quad (a)_n = a(a+1)(a+2) \cdot ... \cdot (a+n-1), \quad n = 1,2,3,....$$

Несколько преобразуем $F_k(z)$ при любых $k$. Меняя порядок суммирования, получаем

$$F_k(z) = \sum_{m=0}^{k} \frac{(-k)_n}{n!} F(-k-1, -k; -k-m; 1) z^m.$$

Пользуясь теперь известной формулой [5, с. 489]

$$F(-n, b; c; 1) = \frac{(c-b)_n}{(c)_n}, \; c-b > -n,$$

имеем

$$F_k(z) = \sum_{m=0}^{k} \frac{(-1)^m (-k)_m (-m)_{k+1}}{m!(k+1)_m} z^m.$$

С учетом формулы $(-m)_{k+1} = 0, \, m \le k$, заключаем, что $F_k(z) \equiv 0$.

Отсюда следует, что оператор $B_{abx}^{n,\lambda}$ – обратный оператору $A_{abx}^{n,\lambda}$, т.е. справедливо первое из равенств теоремы.

Вторая часть теоремы доказывается аналогично.

Имеет место

**Лемма 1.** При $\beta < 1$ и $x \in [a,b]$ справедливы равенства

$$\int_a^x (x-t)^{-\beta} \overline{J}_{-\beta}\left[\lambda\sqrt{(x-t)(b-t)}\right] f(t) dt = \Gamma(1-\beta) D_{ax}^{\beta-1} \left\{ B_{abx}^{1,\lambda}[f(x)] \right\}, \qquad (14)$$

$$\int_x^b (t-x)^{-\beta} \overline{J}_{-\beta}\left[\lambda\sqrt{(t-x)(t-a)}\right] f(t) dt = \Gamma(1-\beta) D_{xb}^{\beta-1} \left\{ B_{bax}^{1,\lambda}[f(x)] \right\}, \qquad (15)$$

где $D_{ax}^l$ и $D_{xb}^l$ – известные операторы дробного интегрирования при $l < 0$ и дробного дифференцирования при $l > 0$; $\overline{J}_\alpha(z) = \Gamma(\alpha+1)(z/2)^{-\alpha} J_\alpha(z)$ – функция Бесселя – Клиффорда.

*Доказательство* равенств (14)-(15) проводится разложением функций Бесселя в степенные ряды и сравнением коэффициентов.

### 3. Обобщенное решение задачи Коши

Решение задачи Коши для уравнения (7), удовлетворяющее начальным условиям (8) и

$$u(\xi, \xi) = \tau(\xi), \; a \le x \le b, \qquad (16)$$

известно [3]:

$$u(\xi, \eta) = \kappa_1 (\eta - \xi)^{-2\beta-1} \int_\xi^\eta (r^\beta \overline{I}_{-\beta}(\lambda\sqrt{r}) - \frac{2\lambda^2}{(1+\beta)(1+2\beta)} r^{1+\beta} \overline{I}_{1+\beta}(\lambda\sqrt{r})) \tau(t) dt -$$

$$- \frac{\kappa_1}{2(1+2\beta)} (\eta - \xi)^{-2\beta-1} \int_\xi^\eta r^\beta \overline{I}_\beta(\lambda\sqrt{r})(\eta + \xi - 2t) \, \tau'(t) dt - \kappa_2 \int_\xi^\eta r^{-\beta} \overline{I}_{-\beta}(\lambda\sqrt{r}) \, \nu(t) dt, \quad (17)$$

где
$$r = |(\eta-t)(t-\xi)|, \quad \kappa_1 = \frac{\Gamma(2+2\beta)}{\Gamma^2(1+\beta)}, \quad \kappa_2 = [2(1-2\beta)]^{2\beta-1}\frac{\Gamma(2-2\beta)}{\Gamma^2(1-\beta)}, \quad \overline{I}_\alpha(z) = \overline{J}_\alpha(iz).$$

Если $\tau(x) \in C^3[a,b]$ и $\nu(x) \in C^2[a,b]$, то функция $u(\xi,\eta)$, определенная формулой (17), является классическим, дважды непрерывно дифференцируемым решением задачи Коши для уравнения (7) с начальными данными (8) и (16) в области $\Delta$.

**Определение 1.** Если функции $\tau'(x)$ и $\nu(x)$ непрерывны при $a < x < b$, то выражение вида (17) будем называть ***обобщенным решением*** уравнения (7) в области $\Delta$.

Для того чтобы обобщенное решение обладало той или иной гладкостью, необходимо, чтобы функции $\tau(x)$ и $\nu(x)$ имели определенную гладкость.

Рассмотрим класс $R_{2k}^\lambda$ обобщенных решений уравнения (7). Здесь и далее $k$ принимает значения $a$ и $b$.

**Определение 2.** ***Обобщенным решением класса*** $R_{2k}^\lambda$ уравнения (7) будем называть функцию $u(\xi,\eta)$ вида (17), где $\tau(x)$ представимо в виде

$$\tau(x) = \tau(k) + \text{sgn}(x-k)\int_k^x |x-t|^{-2\beta}\, \overline{J}_{-\beta}\left[\lambda|x-t|\right]T(t)dt, \tag{18}$$

а $\nu(x)$ и $T(x)$ — непрерывные и интегрируемые в $(a,b)$ функции.

**Замечание 1.** Из (18) нетрудно заключить, что $\tau(x) \in C[a,b]$ и существует $\tau'(x) \in C(a,b)$. Следовательно, обобщенное решение класса $R_{2k}^\lambda$ является обобщенным решением в смысле определения 1.

**Замечание 2.** Не ограничивая общности, будем считать, что $\tau(k) = 0$. При невыполнении этого условия, прибавив к функции $u(\xi,\eta)$ частное решение уравнения (7) вида
$$w(\xi,\eta) = A\, ch(\lambda(\eta+\xi)/2) + B\, sh(\lambda(\eta+\xi)/2),$$
можно распорядиться коэффициентами $A$ и $B$ так, что новая функция в точке $(k,0)$ примет нулевое значение.

Для определенности, обобщенное решение, принадлежащее к $R_{2k}^\lambda$ обозначим через $u_{2k}(\xi,\eta)$.

Согласно [1], интегральные представления обобщенных решений задачи Коши из классов $R_{2a}^\lambda$ и $R_{2b}^\lambda$ соответственно записываются следующим образом:

$$u_{2a}(\xi,\eta) = \int_a^\xi r^{-\beta}\, \overline{J}_{-\beta}(\lambda\sqrt{r})T(t)dt + \int_\xi^\eta r^{-\beta}\, \overline{I}_{-\beta}(\lambda\sqrt{r})N(t)dt; \tag{19}$$

$$u_{2b}(\xi,\eta) = \int_\eta^b r^{-\beta}\, \overline{J}_{-\beta}(\lambda\sqrt{r})T(t)dt + \int_\xi^\eta r^{-\beta}\, \overline{I}_{-\beta}(\lambda\sqrt{r})N(t)dt, \tag{20}$$

где
$$N(t) = [2\cos\beta\pi]^{-1}T(t) - \kappa_2\nu(t).$$

**Лемма 2.** Обобщенное решение $u_{2k}(\xi,\eta) \in R_{2k}^{\lambda}$ обладает следующими свойствами:

1) $u_{2k}(\xi,\eta) \in C(\overline{\Delta}) \cap C^1(\Delta)$, $\dfrac{\partial^2 u}{\partial \xi \partial \eta} \in C(\Delta)$;

2) $u_{2k}(\xi,\eta)$ удовлетворяет уравнению (7) и условиям (8), (16).

Для $\lambda = 0$ лемма 2 доказана в [6]. При наличии $\lambda$ доказательство леммы существенно не отличается.

### 4. Обобщенное решение задачи Коши – Гурса

При решении задачи Коши – Гурса для уравнения (7) будем пользоваться соответствующим интегральным представлением задачи Коши.

Пусть $k = a$. В этом случае, полагая в формуле (19) $\xi = a$ и учитывая условие (9), получаем интегральное уравнение для определения $N(x)$:

$$\int_a^x (x-t)^{-\beta}(t-a)^{-\beta} \overline{I}_{-\beta}\left[\lambda\sqrt{(x-t)(t-a)}\right] N(t) dt = \varphi_a(x).$$

Последнее уравнение в результате применения равенства (14) при $a = b$ можно привести к виду, удобному для дальнейшего исследования:

$$D_{ax}^{\beta-1}\left\{ B_{aax}^{1,\lambda}[(x-a)^{-\beta} N(x)]\right\} = \frac{1}{\Gamma(1-\beta)} \varphi_a(x). \tag{21}$$

Применяя теперь к обеим частям уравнения (21) последовательно операторы $D_{ax}^{1-\beta}$ и $A_{aax}^{1,\lambda}$, получаем

$$T(x) = \kappa_3 \nu(x) + 2\cos\beta\pi \cdot \Phi_a(x), \tag{22}$$

где $\quad \kappa_3 = 2\kappa_2 \cos\beta\pi, \Phi_a(x) = \dfrac{1}{\Gamma(1-\beta)}(x-a)^{\beta} A_{aax}^{1,\lambda}\left\{ D_{ax}^{1-\beta}[\varphi_a(x)]\right\}.$

Подставляя (22) в (19), находим представление решения задачи Коши – Гурса (7) – (9) из класса $R_{2a}^{\lambda}$ в явном виде

$$u_{2a}(\xi,\eta) = \int_a^{\xi} (\xi-t)^{-\beta}(\eta-t)^{-\beta} \overline{J}_{-\beta}\left[\lambda\sqrt{(\eta-t)(\xi-t)}\right][k_3 \nu(t) + 2\cos\beta\pi \cdot \Phi_a(t)] dt +$$

$$+ \int_{\xi}^{\eta} (\eta-t)^{-\beta}(t-\xi)^{-\beta} \overline{I}_{-\beta}\left[\lambda\sqrt{(\eta-t)(t-\xi)}\right] \Phi_a(t) dt. \tag{23}$$

Доказательство, приведенное нами, переносится и на случай $k = b$.

Таким образом, представление решения задачи Коши – Гурса (7), (8), (10) из класса $R_{2b}^{\lambda}$ выписывается в явном виде:

$$u_{2b}(\xi,\eta) = \int_{\eta}^{b} (t-\xi)^{-\beta}(t-\eta)^{-\beta} \overline{J}_{-\beta}\left[\lambda\sqrt{(t-\eta)(t-\xi)}\right][k_3 \nu(t) + 2\cos\beta\pi \cdot \Phi_b(t)] dt +$$

$$+ \int_{\xi}^{\eta} (\eta-t)^{-\beta}(t-\xi)^{-\beta} \overline{I}_{-\beta}\left[\lambda\sqrt{(\eta-t)(t-\xi)}\right] \Phi_b(t) dt, \tag{24}$$

где
$$\Phi_b(x) = \frac{1}{\Gamma(1-\beta)}(b-x)^\beta A_{bbx}^{1,\lambda}\left\{D_{bx}^{1-\beta}[\varphi_b(x)]\right\}.$$

Здесь также можно доказать лемму, аналогичную лемме 2.

**Замечание 3.** Полученные явные интегральные представления (23) и (24) обобщенного решения задачи Коши – Гурса для уравнения (1) играют важную роль при исследовании задач для уравнений смешанного типа, так как из них при $\eta = \xi$ легко вывести основные функциональные соотношения между $\tau(x)$ и $\nu(x)$ на линии вырождения, принесенные из гиперболической части смешанной области.

**Замечание 4.** Класс обобщенных решений уравнения (7) при изучении задачи Коши – Гурса (8), (9) является существенным: если решение уравнения (7), удовлетворяющее условиям (8), (9), не принадлежит к $R_{2a}^\lambda$, то нарушается единственность решения задачи.

**Пример.** Функция [7]
$$u(\xi,\eta) = [(\xi-a)(\eta-a)]^{-\beta}\overline{J}_{-\beta}\left(\lambda\sqrt{(\xi-a)(\eta-a)}\right)$$

является решением уравнения (7), удовлетворяющим однородным условиям (8) и (9). Однако она не принадлежит к классу функций $R_{2a}^\lambda$. В справедливости последнего утверждения можно убедиться, например, с помощью метода от противного.

### 5. Обобщенное решение задачи Гурса

Рассмотрим в области $\Delta$ задачу Гурса для уравнения (7) с условиями (9) – (11). При изучении этой задачи будем пользоваться представлением решения задачи Коши – Гурса для уравнения (7).

Пусть $k = a$. В этом случае, полагая в формуле (23) $\eta = b$ и учитывая условие (10), получаем интегральное уравнение

$$\int_a^x (x-t)^{-\beta}(b-t)^{-\beta}\overline{J}_{-\beta}\left[\lambda\sqrt{(b-t)(x-t)}\right][k_3\nu(t) + 2\cos\beta\pi\cdot\Phi_a(t)]dt +$$
$$+\int_x^b (b-t)^{-\beta}(t-x)^{-\beta}\overline{I}_{-\beta}\left[\lambda\sqrt{(b-t)(t-x)}\right]\Phi_a(t)dt = \varphi_b(x). \quad (25)$$

Разрешая интегральное уравнение (25), имеем
$$k_3\nu(x) + 2\cos\beta\pi\cdot\Phi_a(x) = \frac{(b-x)^\beta}{\Gamma(1-\beta)}A_{abx}^{1,\lambda}\left[D_{ax}^{1-\beta}\varphi_b(x)\right] -$$
$$- (b-x)^\beta A_{abx}^{1,\lambda}\left[D_{ax}^{1-\beta}D_{xb}^{\beta-1}\left\{B_{bbx}^{1,\lambda}[(b-x)^{-\beta}\Phi_a(x)]\right\}\right]. \quad (26)$$

Подставляя (26) в формулу (23), получаем интегральное представление обобщенного решения $u_{2a}(\xi,\eta) \in R_{2a}^\lambda$ задачи Гурса для уравнения (7) с условиями (9) – (11) в виде

$$u_{2a}(\xi,\eta) = \frac{1}{\Gamma(1-\beta)}\int_a^\xi (\xi-t)^{-\beta}(\eta-t)^{-\beta}\overline{J}_{-\beta}\left[\lambda\sqrt{(\eta-t)(\xi-t)}\right](b-t)^\beta A_{abt}^{1,\lambda}\left\{D_{at}^{1-\beta}[\varphi_b(t)]\right\}dt -$$

$$-\int_a^\xi (\xi-t)^{-\beta}(\eta-t)^{-\beta}\overline{J}_{-\beta}\left[\lambda\sqrt{(\eta-t)(\xi-t)}\right](b-t)^\beta A_{abt}^{1,\lambda}\left\{D_{at}^{1-\beta}\left[D_{tb}^{\beta-1}\left\{B_{bbt}^{1,\lambda}[(b-t)^{-\beta}\Phi_a(t)]\right\}\right]\right\}dt +$$

$$+\int_\xi^\eta (\eta-t)^{-\beta}(t-\xi)^{-\beta}\overline{I}_{-\beta}\left[\lambda\sqrt{(\eta-t)(t-\xi)}\right]\Phi_a(t)dt. \qquad (27)$$

Аналогично находится обобщенное решение этой же задачи, принадлежащее к $R_{2b}^\lambda$:

$$u_{2b}(\xi,\eta)=\frac{1}{\Gamma(1-\beta)}\int_\eta^b (t-\xi)^{-\beta}(t-\eta)^{-\beta}\overline{J}_{-\beta}\left[\lambda\sqrt{(t-\eta)(t-\xi)}\right](t-a)^\beta A_{bat}^{1,\lambda}\left\{D_{tb}^{1-\beta}[\varphi_a(t)]\right\}dt -$$

$$-\int_\eta^b (t-\xi)^{-\beta}(t-\eta)^{-\beta}\overline{J}_{-\beta}\left[\lambda\sqrt{(t-\eta)(t-\xi)}\right](t-a)^\beta A_{bat}^{1,\lambda}\left\{D_{tb}^{1-\beta}\left[D_{at}^{\beta-1}\left\{B_{aat}^{1,\lambda}[(t-a)^{-\beta}\Phi_b(t)]\right\}\right]\right\}dt +$$

$$+\int_\xi^\eta (\eta-t)^{-\beta}(t-\xi)^{-\beta}\overline{I}_{-\beta}\left[\lambda\sqrt{(\eta-t)(t-\xi)}\right]\Phi_b(t)dt. \qquad (28)$$

**Замечание 5.** Из формул (27) и (28) нетрудно заметить, что при исследовании задачи Гурса для уравнений гиперболического типа второго рода, в отличии от уравнений первого рода, нарушается равноправие характеристик. Это обстоятельство связано с необходимостью введения представления (18).

In this paper, the Cauchy, Cauchy–Goursat, and Goursat problems for a degenerate second kind hyperbolic equation with a spectral parameter are studied. For these equations, depending on the degree of degeneracy, limit values of the sought solutions and its derivative on degeneration lines can have singularities. To provide the required smoothness of the solution outside the characteristic line of degeneration, it is necessary to require enhanced data smoothness. In order to ease this requirement, a definition of a class of generalized solutions is introduced and properties of this class are studied. In addition, on the basis of the well-known formula of the classical solution of the Cauchy problem for the above equation, a generalized solution of the Cauchy problem in the introduced class is obtained in an explicit form which is easy to use for further research. Properties of these solutions are studied. Some operators with Bessel functions in the nucleus are introduced and their basic properties are studied. The proved important identities of these operators and the above representation of the generalized solution of the Cauchy problem allow one to find an explicit representation of the generalized solutions of the Cauchy–Goursat and Goursat problems in the characteristic triangle. In addition, an example showing the importance of introducing such class is presented: if the solution does not belong to the newly introduced class, then the uniqueness of the solution of the Cauchy–Goursat problem can be broken. The resulting explicit integral representation of the generalized solution of the Cauchy–Goursat problem plays an important role in the study of problems for equations of the mixed type: it makes it easy to derive the basic functional relationship between the traces of the sought solution and of its derivative on the line of degeneration from the hyperbolic part of the mixed domain.

Keywords: degenerate hyperbolic equation of the second kind, the spectral parameter, generalized solution, the operator with the Bessel functions in the nucleus.

*EHRGASHEV Tuhtasin Gulamzhanovich* (Tashkent Institute of Irrigation and Melioration (TIIM), Tashkent, Uzbekistan)
E-mail: ertuhtasin@mail.ru